\documentclass[11pt]{amsart}
\usepackage{amsmath, amssymb, amsthm, bm, mathtools, enumitem, caption}
\usepackage{hyperref}

\usepackage[noabbrev,capitalize]{cleveref}
\usepackage{adjustbox}
\crefname{equation}{}{}
\usepackage{ytableau}
\usepackage{graphicx, tikz}
\usepackage[textheight=7.6in, textwidth=6.5in]{geometry}

\usepackage{microtype}

\newtheorem{theorem}{Theorem}[section]
\newtheorem{lemma}[theorem]{Lemma}
\newtheorem{corollary}[theorem]{Corollary}
\newtheorem{proposition}[theorem]{Proposition}

\newtheorem*{conjecture*}{Conjecture}

\theoremstyle{definition}

\theoremstyle{remark}
\newtheorem*{remark}{Remark}

\newtheorem*{example}{Example}

\numberwithin{equation}{section}

\DeclareMathOperator{\Tr}{Tr}

\newcommand{\fmodd}{\mathfrak{m}_{{\text {\rm odd}}}}

\title[Recursive Formulas]{Recursive Formulas for MacMahon and Ramanujan $q$-series}
\author{Tewodros Amdeberhan, Rupam Barman \and Ajit Singh}
\address{Dept. of Mathematics, Tulane University, New Orleans, LA 70118, USA}
\email{tamdeber@tulane.edu}
\address{Department of Mathematics, Indian Institute of Technology Guwahati, Assam, India \\  \,\,\,\,  PIN- 781039}
\email{rupam@iitg.ac.in}
\address{Dept. of Mathematics, University of Virginia, Charlottesville, VA 22904, USA}
\email{ajit18@iitg.ac.in}

\begin{document}
\begin{abstract} 
In the present work, we extend current research in a nearly-forgotten but newly revived topic, initiated by P. A. MacMahon, on a generalized notion which relates the divisor sums to the theory of integer partitions and two infinite families of $q$-series by Ramanujan. Our main emphasis will be on explicit representations for a variety of $q$-series, studied primarily by MacMahon and Ramanujan, with an eye towards their modular properties and their proper place in the ring of quasimodular forms of level one and level two.
 \end{abstract}

\dedicatory{Dedicated to George Andrews and Bruce Berndt on their $1010101$-th birthdays}

\maketitle
\section{Introduction and Statement of Results}

\noindent
The classical sequence of Eisenstein series are defined as (for example, see Chapter 1 of \cite{CBMS}) 
\begin{equation}\label{Eisenstein}
E_{2k}(q):=1-\frac{4k}{B_{2k}}\sum_{n=1}^{\infty}\sigma_{2k-1}(n)q^n,
\end{equation}
where $B_{k}$ is the $k$-th Bernoulli number and $\sigma_k(n)=\sum_{d\vert n}d^k$ is the power-sum of divisors of $n$.

\smallskip
\noindent
In the present work, we chose to follow up on Ramanujan's predictions regarding two functions \cite[page 369]{Rama} that he himself defined. Namely that
\begin{align*}
	U_{2t}(q)&=\frac{1^{2t+1}-3^{2t+1}q+5^{2t+1}q^3-7^{2t+1}q^6+9^{2n+1}q^{10}-\cdots}
	{1-3q+5q^3-7q^6+9q^{10}-\cdots} \qquad \text{and} \\
	V_{2t}(q)&=\frac{1^{2t}-5^{2t}q-7^{2t}q^2+11^{2t}q^5+13^{2t}q^7-\cdots}
	{1-q-q^2+q^5+q^7-\cdots}.
\end{align*}
After listing the first few expansions $U_0=V_0=1, U_2=V_2=E_2$,
\begin{align*}
	U_4&=\frac13(5E_2^2-2E_4), \qquad \qquad \qquad \qquad V_4=3E_2^2-2E_4, \\
	U_6&=\frac19(35E_2^3-42E_2E_4+16E_6), \qquad \,\,\, V_6=15E_2^2-30E_2E_4+16E_6, \\
	U_8&=\frac13(35E_2^4-84E_2^2E_4+64E_2E_6-12E_4^2), \,\,\,\,
	V_8=105E_2^4-420E_2^2E_4+448E_2E_6-132E_4^2. 
\end{align*}
Ramanujan declared that: 
$$\text{``\emph{In general $U_{2t}$ and $V_{2t}$ are of the form
		$\sum K_{\ell,m,n} \, E_2^{\ell}E_4^mE_6^n,$ where $\ell+2m+3n=t$."}}$$ 
In modern language Ramanujan's declaration amounts to saying that both $U_{2t}$ and $V_{2t}$  are quasimodular forms of weight $2t$. 
Berndt, Chan, Liu, Yee, and Yesilyurt \cite{Berndt1, Berndt2} proved this claim using Ramanujan's identities \cite{Rama2}
\begin{equation}  \label{diffeq}
D(E_2) =\frac{E_2^2-E_4}{12},\ \  \
D(E_4)=\frac{E_2E_4-E_6}{3},\ \ \ {\text {\rm and}}\ \ \
D(E_6)=\frac{E_2E_6-E_4^2}{2},
\end{equation}
where $D:=q\frac{d}{dq}.$ 
However, their results are not explicit. Indeed,
Andrews and Berndt (see p. 364 of \cite{AndrewsBerndt}) proclaim that \emph{``...it seems extremely difficult to find a general formula for all $K_{\ell,m,n}.$''} 

\smallskip
\noindent
Ramanujan's claim is that $U_{2t}(q)$ and $V_{2t}(q)$  are weight $2t$ quasimodular forms, as  the ring of quasimodular forms  is the polynomial ring  (for example, see \cite{Zagier})
$$
\mathbb{C}[E_2,E_4,E_6]=\mathbb{C}[E_2, E_4, E_6, E_8, E_{10}, \dots].
$$
The first and the third authors together with K. Ono \cite{AOS2} chose the latter $\mathbb{C}[E_2, E_4, E_6, E_8, E_{10},\dots]$, involving all Eisenstein series expressed as ``traces of partition Eisenstein series,'' to produce the first explicit formulas for both $U_{2t}(q)$ and $V_{2t}(q)$.

\smallskip
\noindent
Our goal in this paper is to obtain recursive formulas in the ring $\mathbb{C}[E_2,E_4,E_6]$ just as Ramanujan originally proposed.
One can argue that the main merit of our effort here lies in inviting the audience to a variety to the techniques employed for the present goal, that the authors believe 
should help in similar circumstances. The first of such installments appeared in \cite{AOS} for the $q$-series $U_{2t}(q)$:
\begin{theorem} \label{theorem1}
	If $t$ is a non-negative integer, then we have that
	$$
	U_{2t}(q)=
	\sum_{\substack {\alpha, \beta, \gamma\geq 0\\ \alpha+2\beta+3\gamma=t}}
	c_u(\alpha,\beta,\gamma) E_2(q)^{\alpha} E_4(q)^{\beta} E_6(q)^{\gamma}
	$$
where the coefficients $c_u(\alpha,\beta,\gamma)$ are rational numbers defined by \cite[eq'n (1.7)]{AOS}.
\end{theorem}

\noindent
Our first result at present concerns the other Ramanujan function $V_{2t}(q)$. We require in introductions of a triple-indexed sequence of \emph{integers} defined by 
\begin{align} \label{coefficientsV} \nonumber
	c_v(\alpha,\beta,\gamma)
	&:=(2\alpha+8\beta+12\gamma-1)\cdot c_v(\alpha-1,\beta,\gamma) - 2(\alpha+1)\cdot c_v(\alpha+1,\beta-1,\gamma) \\
	&\qquad  \ \ \ \ \  - 8(\beta+1)\cdot c_v(\alpha,\beta+1,\gamma-1)- 12(\gamma+1)\cdot c_v(\alpha,\beta-2,\gamma+1),
\end{align}
where $\alpha, \beta, \gamma\geq0.$ To seed the recursion, we let $c_v(0, 0, 0):=1,$ and we let
$c_v(\alpha,\beta,\gamma):=0$ if any of the arguments are negative. 
Here we list the ``first few'' values: 
\begin{align*}
	c_v(1,0,0)=1, \  c_v(0,1,0)=-2,\  c_v(0,0,1)=16,  \   c_v(1,1,0)=-30,  \  
	c_v(1,0,1)=448,\dots.
\end{align*}

\begin{theorem}\label{theorem2}
	If $t$ is a non-negative integer, then we have that
	$$
	V_{2t}(q)=
	\sum_{\substack {\alpha, \beta, \gamma\geq 0\\ \alpha+2\beta+3\gamma=t}}
	c_v(\alpha,\beta,\gamma) E_2(q)^{\alpha} E_4(q)^{\beta} E_6(q)^{\gamma}
	$$
where the coefficients $c_v(\alpha,\beta,\gamma)$ are defined by \eqref{coefficientsV}.
\end{theorem}

\smallskip
\noindent
In the same spirit but in an earlier paper \cite{AOS}, the first and the third authors together with K. Ono have found such an explicit description for MacMahon's quasimodular form 
\begin{align} \label{MO_U}
\mathcal{U}_{2t}(q):=\sum_{1\leq k_1<\cdots<k_t} \frac{q^{k_1+\cdots+k_t}}{(1-q^{k_1})^2\cdots(1-q^{k_t})^2}.
\end{align}
An important observation is that Ramanujan's $U_{2t}(q)$ and MacMahon's $\mathcal{U}_{2t}(q)$ are directly linked to each other. In fact, this fact allows us to formulate Theorem 1.3 of \cite{AOS} as follows:
$$\mathcal{U}_{2t}(q)=\sum_{a=0}^t w_a(t)\cdot U_{2a}(q) \qquad \text{where} \qquad 
w_a(t):=\frac{\binom{2t}t}{16^t(2t+1)}\sum_{0\leq\ell_1<\cdots<\ell_a<t} \prod_{j=1}^a\frac1{(2\ell_j+1)^2}.$$
In light of this, the weight $2a$ part of $\mathcal{U}_{2t}(q)$ becomes precisely $U_{2a}(q)$ and we can evidently borrow the corresponding expansion from Theorem~\ref{theorem1} above.

\smallskip
\noindent
In a related rendition, the first author together with Andrews and Tauraso \cite{AAT} introduced a $q$-series which is intimately connected to MacMahon's $\mathcal{U}_{2t}(q)$ and given by
\begin{align} \label{AAT_U}
\mathcal{U}_{2t}^{\star}(q):=\sum_{1\leq k_1\leq \cdots\leq k_t} \frac{q^{k_1+\cdots+k_t}}{(1-q^{k_1})^2\cdots(1-q^{k_t})^2}.
\end{align}
In the same paper \cite[Theorem 6.1]{AAT}, the authors have shown that each of these $\mathcal{U}_{2t}^{\star}(q)$ are quasimodular forms of weight at most $2t$. 
The first and the third authors together with K. Ono followed this through and furnished the expansion \cite[Theorem 1.4]{AOS}
$$\mathcal{U}_{2t}^{\star}(q)=\sum_{a=0}^t w_a^{\star}(t)\cdot \mathbb{E}_{2a}^{\star}(q) \qquad \text{where} \qquad
\mathbb{E}_{2a}^{\star}(q):=\sum_{(1^{m_1},\dots,a^{m_a})\vdash a}\,
\prod_{j=1}^a\frac1{m_j!}\left(-\frac{B_{2j}\,E_{2j}(q)}{(2j)\cdot (2j!)}\right)^{m_j}$$
and for some coefficients $w_a^{\star}(t)$ akin to the above $w_a(t)$.

\smallskip
\noindent 
In the interest of exhibiting the tight link between the two $q$-series, $\mathcal{U}_{2t}(q)$ and $\mathcal{U}_{2t}^{\star}(q)$, we brought to bear the following relation which expresses one in terms of the other:
\begin{align*}
\mathcal{U}_{2t}^{\star}(q)&=(-1)^t\sum_{(1^{m_1},\cdots, t^{m_t})\vdash t}
(-1)^{m_1+\cdots+m_t}\binom{m_1+\cdots+m_t}{m_1,\dots,m_t}
\prod_{k=1}^t\left(\mathcal{U}_{2k}(q)\right)^{m_k}.
\end{align*}
The proof relies on a well-known convolution between the elementary symmetric functions $\mathbf{e}_t$ and the complete homogeneous symmetric functions $\mathbf{h}_t$,  effectively  utilized in \cite[Lemma 6.1]{AAT}:
$$\sum_{i=0}^t(-1)^i \mathbf{e}_i\mathbf{h}_{t-i}=0.$$

\smallskip
\noindent
In addition to the $\mathcal{U}_{2t}(q),$ MacMahon also introduced \cite{MacMahon}  the $q$-series
\begin{align} \label{Mac2} 
	\mathcal{C}_{2t}(q)=\sum_{n=1}^{\infty} \fmodd(t;n)q^n:=\sum_{0< s_1<s_2<\cdots<s_t} \frac{q^{2s_1+2s_2+\cdots+2s_t-t}}{(1-q^{2s_1-1})^2(1-q^{2s_2-1})^2\cdots(1-q^{2s_t-1})^2}.
\end{align}
The number $\fmodd(t;n)$ is a sum of the products of the {\it part multiplicities}  for partitions of $n$ with $t$ distinct odd part sizes.  
Furthermore, in analogy with the work of Andrews and Rose \cite{Andrews-Rose, Rose},  Bachmann \cite{Bachmann} proved that each $\mathcal{C}_{2t}(q)$ is a finite linear combination of quasimodular forms on $\Gamma_0(2)$ of weight at most $2t$. In recent years, the literature saw a flurry of research activities related to the present context, see 
\cite{AAT2, AGOS, JPS, Ono-Singh}.

\smallskip
\noindent
Here we give yet another explicit formula for $\mathcal{C}_{2t}(q)$  analogous to the MacMahon's 
$\mathcal{U}_{2t}(q)$ series as described in \cite[Theorem 1.3]{AOS}. We require some preliminary concepts and terminologies.
Let's recall the Jacobi theta series
$$\theta_2(q)=\sum_{m\in\mathbb{Z}+\frac12}q^{m^2}, \qquad \theta_3(q)=\sum_{m\in\mathbb{Z}}q^{m^2} \qquad \text{and} \qquad
\theta_4(q)=\sum_{m\in\mathbb{Z}} (-1)^mq^{m^2},$$
satisfying the identity $\theta_3^4(q)=\theta_2^4(q)+\theta_4^4(q)$ and having the derivatives
\begin{align} \label{Ram_like}
\frac{D(\theta_2)}{\theta_2}&=\frac{E_2-\theta_2^4+5\theta_3^4}{24}, \qquad
\frac{D(\theta_3)}{\theta_3}=\frac{E_2+5\theta_2^4-\theta_3^4}{24}, \qquad
\frac{D(\theta_4)}{\theta_4}=\frac{E_2(q) - \theta_2^4-\theta_3^4}{24}.
\end{align}
Observe that modular forms over the congruence subgroup $\Gamma_0(2)$ can be generated by $\Theta_{0,1}(q)$ and $\Theta_{1,1}(q)$ where these symmetric combinations are defined by
$$\Theta_{r,s}(q):=\theta_2^{4r}(q)\cdot \theta_3^{4s}(q)+\theta_2^{4s}(q)\cdot \theta_3^{4r}(q)$$
which carry weight $2r+2s$.
Construct the sequences $c_c(\alpha,\beta,\gamma)$ defined recursively by
\begin{align} \label{coeff_C} 
c_c(\alpha,\beta,\gamma) \nonumber
&=\frac{2\alpha+4\beta+4\gamma-1}{24}c_c(\alpha-1,\beta,\gamma)+\frac{20\gamma-4\beta+3}{24}c_c(\alpha,\beta-1,\gamma) \\   
& + \frac{20\beta-4\gamma+3}{24}c_c(\alpha,\beta,\gamma-1) - \frac{7(\alpha+1)}6 c_c(\alpha+1,\beta-1,\gamma-1) \\   \nonumber
& - \frac{\alpha+1}{12}c_c(\alpha+1,\beta-2,\gamma) - \frac{\alpha+1}{12}c_c(\alpha+1,\beta,\gamma-2). 
\end{align}
To seed the recursion, we let $c_c(\alpha,\beta,\gamma) := 0$ if any of the arguments is negative while we let $c_c(0, 0, 0) := 1$.
We are now ready to state the content of the promised results.

\begin{theorem} \label{theorem3} If $t$ is a positive integer, then we have
$$\frac{D^t(\theta_4(q))}{\theta_4(q)} 
=\sum_{\substack{\alpha, \beta, \gamma\geq0 \\ 2\alpha+2\beta+2\gamma=2t}}
c_c(\alpha,\beta,\gamma)\cdot E_2^{\alpha}(q)\,\Theta_{\beta,\gamma}(q)$$
where the coefficients $c_c(\alpha,\beta,\gamma)$ are defined by \eqref{coeff_C}.
\end{theorem}

\begin{example}  In view of Theorem~\ref{theorem3}, we calculate the following weight $4$ modular form as:
\begin{align} 
\frac{D^2(\theta_4)}{\theta_4}&=\frac1{192}E_2^2-\frac1{96}E_2\Theta_{0,1}+\frac1{192}\Theta_{0,2}-\frac{11}{96}\Theta_{1,1}.
\end{align}
\end{example}

\smallskip
\noindent
We require the role of the \emph{elementary symmetric functions} given by $\mathbf{e}_0=1$ and 
$$\mathbf{e}_k(x_1,\dots,x_r)=\sum_{1\leq i_1<\cdots < i_k\leq r}x_{i_1}\cdots x_{i_k} \qquad \text{for $1\leq k\leq r$}.$$

\begin{theorem} \label{corollary1}
If $t$ is a positive integer, then we have that
$$\mathcal{C}_{2t}(q)=\sum_{k=1}^t (-1)^k\,\frac{\mathbf{e}_{t-k}(0^2,1^2,\dots,(t-1)^2)}{(2t)!}
\sum_{\substack{\alpha, \beta, \gamma\geq0 \\ \alpha+\beta+\gamma=t}}
c_c(\alpha,\beta,\gamma)\cdot E_2^{\alpha}(q)\,\Theta_{\beta,\gamma}(q).$$
\end{theorem}

\smallskip
\noindent
This paper is organized as follows. In Section 2, we prove Theorem \ref{theorem2}. In Section 3, we prove Theorem~\ref{theorem3} and by using the theory of elementary symmetric function, we prove Theorem \ref{corollary1}. In Section 4, we give the generalization of Ramanujan's derivatives \eqref{diffeq} and in Section 5, we give sketch of another proof of Theorem \ref{theorem2}.

\section*{Acknowledgements}
The third  author is grateful for the support of a Fulbright Nehru Postdoctoral Fellowship.

\section{Proof of Theorem  \ref{theorem2}}

\begin{proof}[Proof of Theorem \ref{theorem2}] 
Let $\eta(q):=q^{\frac1{24}}\prod_{k\geq1}(1-q^k)$ be the Dedekind's eta function. We calculate $V_{2t}(q)=\frac{D^t(\eta(q))}{\eta(q)}$ by inducting on $t$. 
	First observe $D(\eta(q))=\frac1{24}\eta(q) E_2(q)$. The functions $V_{2t}(q)$ are known  (see, for example \cite{AOS2, Berndt1, Berndt2}) to be quasimodular forms of weight $2t$, over $\text{SL}_2(\mathbb{Z})$, thus we are ensured that there exist numbers $\widetilde{c}_v(\alpha,\beta,\gamma)$ for which
	$$\frac{D^t(\eta(q))}{\eta(q)}=\sum_{\substack {\alpha, \beta, \gamma\geq0 \\ \alpha+2\beta+3\gamma=t}} 
	\widetilde{c}_v(\alpha,\beta,\gamma)
	\cdot E_2^{\alpha}(q)E_4^{\beta}(q)E_6^{\gamma}(q).$$
     Writing $\widetilde{c}_v$ instead of $\widetilde{c}_v(\alpha,\beta,\gamma)$, one more derivative $D=q\frac{d}{dq}$ turns the last equation into 
	\begin{align*}
		D^{t+1}(\eta)
		&=D (\eta)\cdot\left(\sum_{\alpha,\beta,\gamma}\widetilde{c}_v \cdot E_2^{\alpha}E_4^{\beta}E_6^{\gamma}\right)
		+\eta\cdot \sum_{\alpha, \beta, \gamma}\widetilde{c}_v\cdot D(E_2^{\alpha}E_4^{\beta}E_6^{\gamma}).
	\end{align*}
	On the other hand, Ramanujan's identities \eqref{diffeq} imply that
	\begin{align*}
		D(E_2^{\alpha}E_4^{\beta}E_6^{\gamma})
		&=\left(\frac{\alpha}{12}+\frac{\beta}3+\frac{\gamma}2\right)E_2^{\alpha+1}E_4^{\beta}E_6^{\gamma}-\frac{\alpha}{12}E_2^{\alpha-1}E_4^{\beta+1}E_6^{\gamma}  \\
		& \qquad \qquad \qquad - \frac{\beta}3E_2^{\alpha}E_4^{\beta-1}E_6^{\gamma+1}  - \frac{\gamma}2E_2^{\alpha}E_4^{\beta+2}E_6^{\gamma-1}.
	\end{align*}
	We find that the homogeneous weight $2t+2$ form $D^{t+1}(\eta)$ satisfies
	\begin{align*}
		\frac{D^{t+1}(\eta)}{\eta}
		&=\sum_{\substack {\alpha, \beta, \gamma\geq0 \\ \alpha+2\beta+3\gamma=t}} 
		\left(\frac{\alpha}{12}+\frac{\beta}3+\frac{\gamma}2+\frac1{24}\right)
		\widetilde{c}_v\cdot E_2^{\alpha+1}E_4^{\beta}E_6^{\gamma} 
		-\sum_{\alpha,\beta,\gamma}\frac{\alpha}{12} \,\widetilde{c}_v\cdot E_2^{\alpha-1}E_4^{\beta+1}E_6^{\gamma}  \\
		& \qquad \qquad \qquad - \sum_{\alpha,\beta,\gamma}\frac{\beta}3 \, \widetilde{c}_v\cdot E_2^{\alpha}E_4^{\beta-1}E_6^{\gamma+1}  
		- \sum_{\alpha,\beta,\gamma}\frac{\gamma}2 \, \widetilde{c}_v\cdot E_2^{\alpha}E_4^{\beta+2}E_6^{\gamma-1}.
	\end{align*}
	By comparing the coefficients of $E_2^{\alpha}E_4^{\beta}E_6^{\gamma}$ on both sides of the equation above,
	we obtain the recursion  (together with $\widetilde{c}_v(\alpha,\beta,\gamma)=\delta_{(0,0,0)}(\alpha,\beta,\gamma)$, a Kronecker delta)
	\begin{align*}
		\widetilde{c}_v(\alpha,\beta,\gamma)
		&=\left(\frac{\alpha}{12}+\frac{\beta}3+\frac{\gamma}2-\frac1{24}\right)\widetilde{c}_v(\alpha-1,\beta,\gamma)
		-\frac{\alpha+1}{12}\cdot \widetilde{c}_v(\alpha+1,\beta-1,\gamma) \\
		&\qquad \qquad \qquad -\frac{\beta+1}3\cdot \widetilde{c}_v(\alpha,\beta+1,\gamma-1) - \frac{\gamma+1}2\cdot 
		\widetilde{c}_v(\alpha,\beta-2,\gamma+1).
	\end{align*}
To determine the exact weight $2t$ term, we take into account the
	prefactor $24^{\alpha+2\beta+3\gamma}$ so that 
$$V_{2t}(q)=24^{2t}\frac{D^t(\eta(q))}{\eta(q)} \qquad \text{and} \qquad 
c_v(\alpha,\beta,\gamma):=24^{\alpha+2\beta+3\gamma}\cdot \widetilde{c}_v(\alpha,\beta,\gamma).$$
As a result, we obtain the desired
	\begin{align*}
		c_v(\alpha,\beta,\gamma)
		&=\left(2\alpha+8\beta+12\gamma-1\right)\cdot c_v(\alpha-1,\beta,\gamma)
		- 2(\alpha+1) \cdot c_v(\alpha+1,\beta-1,\gamma) \\
		&\qquad \qquad - 8(\beta+1) \cdot c_v(\alpha,\beta+1,\gamma-1) - 12(\gamma+1) \cdot c_v(\alpha,\beta-2,\gamma+1).
	\end{align*}
	We thus conclude the constructions and the proof.
\end{proof}

\section{Proof of Theorems \ref{theorem3} and \ref{corollary1}}

\begin{proof}[Proof of Theorem \ref{theorem3}]
We proceed by induction on $t\geq1$. Since $\Theta_{0,1}=\theta_2^4+\theta_3^4$, the base case $t=1$ is recovered from \eqref{Ram_like} with $c_c(1,0,0)=\frac1{24}$ and 
$c_c(0,0,1)=-\frac1{24}$. Also, note the additional properties that
\begin{align} \label{prelim1}
\Theta_{a,b}\Theta_{a',b'}&=\Theta_{a+a',b+b'}+\Theta_{a+b',b+a'}, \qquad \Theta_{a,b}=\Theta_{b,a}, \qquad E_4=\Theta_{0,2}+7\Theta_{1,1}.
\end{align}
Moreover, the following relation follows from \eqref{Ram_like}: 
\begin{align} \label{prelim2}
D(\Theta_{r,s})&=\frac{r+s}6E_2\Theta_{r,s}+\frac{5s-r}6\Theta_{r+1,s}+\frac{5r-s}6\Theta_{r,s+1}.
\end{align}
Assume the expansion of the weight $2t$ quasimodular form
$$\frac{D^t(\theta_4)}{\theta_4}=\sum_{\substack{\alpha, \beta, \gamma\geq0 \\ \alpha+\beta+\gamma=t}} c_c(\alpha,\beta,\gamma)\cdot E_2^{\alpha}\Theta_{\beta,\gamma}.$$
Then, we take one more derivative and apply \eqref{Ram_like}, \eqref{prelim1}, \eqref{prelim2} to obtain
\begin{align*} 
\frac{D^{t+1}(\theta_4)}{\theta_4}&=\frac{D(\theta_4)}{\theta_4}\left(\sum c_c\cdot E_2^{\alpha}\Theta_{\beta,\gamma}\right)
+\sum c_c\cdot [\Theta_{\beta,\gamma} D(E_2^{\alpha})+E_2^{\alpha}D(\Theta_{\beta,\gamma})] \\
&=\left(\frac{E_2(q)-\Theta_{0,1}}{24}\right)\left(\sum c_c\cdot E_2^{\alpha}\Theta_{\beta,\gamma}\right) 
+ \sum c_c\cdot \left[\frac{\alpha \Theta_{\beta,\gamma} E_2^{\alpha-1}(E_2^2-\Theta_{0,2}-7\Theta_{1,1})}{12}  \right] \\
& + \sum c_c\cdot \left[\frac{E_2^{\alpha}\cdot [(\beta+\gamma)E_2\Theta_{\beta,\gamma}+(5\gamma-\beta)\Theta_{\beta+1,\gamma}+(5\beta-\gamma)
\Theta_{\beta,\gamma+1}]}6 \right]  \\
&=\frac1{24}\sum c_c\cdot [ E_2^{\alpha+1}\Theta_{\beta,\gamma} - E_2^{\alpha}\Theta_{\beta+1,\gamma}-E_2^{\alpha}\Theta_{\beta,\gamma+1}]\\
& + \frac1{12}\sum c_c\cdot \left[\alpha E_2^{\alpha+1} \Theta_{\beta,\gamma}  
- \alpha E_2^{\alpha-1}(\Theta_{\beta+2,\gamma}+\Theta_{\beta,\gamma+2}) - 14\alpha E_2^{\alpha-1}\Theta_{\beta+1,\gamma+1} \right] \\
& + \frac16 \sum c_c\cdot \left[(\beta+\gamma) E_2^{\alpha+1}\Theta_{\beta,\gamma}+(5\gamma-\beta)E_2^{\alpha}\Theta_{\beta+1,\gamma}+(5\beta-\gamma)
E_2^{\alpha}\Theta_{\beta,\gamma+1}] \right] \\
&=\frac1{24}\sum (2\alpha+4\beta+4\gamma+1)\cdot c_c\cdot E_2^{\alpha+1}\Theta_{\beta,\gamma} \\
& + \frac1{24} \sum (20\gamma-4\beta-1)\cdot c_c\cdot E_2^{\alpha}\Theta_{\beta+1,\gamma} \\
& + \frac1{24} \sum (20\beta-4\gamma-1)\cdot c_c\cdot E_2^{\alpha}\Theta_{\beta,\gamma+1}  \\
& - \frac1{12}\sum \alpha \cdot c_c\cdot \left[
 E_2^{\alpha-1} \Theta_{\beta+2,\gamma}+ E_2^{\alpha-1} \Theta_{\beta,\gamma+2} + 14 E_2^{\alpha-1}\Theta_{\beta+1,\gamma+1} \right];
\end{align*}
where we wrote $c_c$ for $c_c(\alpha,\beta,\gamma)$ assuming no confusion arises.
By comparing coefficients on both sides of the last equation, one is lead to the recurrence
\begin{align*}
c_c(\alpha,\beta,\gamma)
&=\frac{2\alpha+4\beta+4\gamma-1}{24}c_c(\alpha-1,\beta,\gamma)+\frac{20\gamma-4\beta+3}{24}c_c(\alpha,\beta-1,\gamma) \\
& + \frac{20\beta-4\gamma+3}{24}c_c(\alpha,\beta,\gamma-1) - \frac{7(\alpha+1)}6 c_c(\alpha+1,\beta-1,\gamma-1) \\
& - \frac{\alpha+1}{12}c_c(\alpha+1,\beta-2,\gamma) - \frac{\alpha+1}{12}c_c(\alpha+1,\beta,\gamma-2). 
\end{align*}
The proof is hence complete.
\end{proof}

\begin{proof}[Proof of Theoerem \ref{corollary1}]
Begin by expressing the quasimodular form $\mathcal{C}_{2t}(q)$ in the manner
$$\mathcal{C}_{2t}(q)=\sum_{k=0}^t (-1)^k \frac{v_t(k)}{(2t)!}\cdot \frac{D^k(\theta_4(q)}{\theta_4(q)}.$$
A direct utility of the relation (see \cite[Corollary 3]{Andrews-Rose})
$$\mathcal{C}_{2t}(q)=\frac1{2t(2t-1)}\left[ (2\mathcal{C}_1(q)+(t-1)^2)\mathcal{C}_{2t-2}(q)-D\mathcal{C}_{2t-2}(q)\right],$$
implies the recurrence  $v_t(k)=(t-1)^2v_{t-1}(k)+v_{t-1}(k-1)$. However, one easily checks that
$$v_t(k)=\mathbf{e}_{t-k}(0^2,1^2,\dots,(t-1)^2)$$
is indeed the solution. The proof follows from Theorem~\ref{theorem3}.
\end{proof}

\bigskip

\section{Generalizing Ramanujan's derivatives for $E_2, E_4$ and $E_6$}

\noindent
Ramanujan famously obtained the following formulas \cite[p. 181]{Rama2} for the action of $D=q\frac{d}{dq}$:
\begin{align*}  
	D (E_2)&=\frac{E_2^2-E_4}{12}, \qquad D(E_4)=\frac{E_2E_4-E_6}3,  \qquad
	D (E_6)=\frac{E_2E_6-E_4^2}2. 
\end{align*}

\noindent
Let $B_s$ denote the Bernoulli numbers and recall the partition Eisenstein series \cite[eq'n 1.5]{AOS2}
\begin{align*}
	\mathbb{E}_{2t}(q)&:=\sum_{(1^{m_1},\dots,t^{m_t})\vdash t}
	\prod_{s=1}^t\frac1{m_s!}\left(\frac{B_{2s}\,E_{2s}(q)}{(2s)\cdot (2s)!}\right)^{m_s}.
\end{align*}
These series play a valuable role (see, for example \cite[Theorem 1.4 and eq'n (3.8)]{AOS}) in the identification of weight $2t$ components of the quasimodular forms $\mathcal{U}_{2a}$ and $\mathcal{U}_{2a}^{\star}$, and they are \emph{universal} (i.e. independent of the ``$a$'' up to a constant factor).

\smallskip
\noindent
Remember the Dedekind's $\eta$-function can be given by $\eta(q)=q^{\frac1{24}}\prod_{k\geq1}(1-q^k)$. Now, denote
$\psi(q)=\eta^3(q)$ and $\mathcal{E}_t(q)=8^t\frac{D^t(\psi(q))}{\psi(q)}$.
In particular, $\mathcal{E}_0=1$ and $\mathcal{E}_1=E_2$.

\begin{lemma} \label{diff-for-F} The functions $\mathcal{E}_t(q)$ satisfy the differential equation 
	$\mathcal{E}_t(q)=(E_2(q)+8D)\mathcal{E}_{t-1}(q)$. 
\end{lemma}
\begin{proof} Employing $8\frac{D(\psi)}{\psi}=E_2=\mathcal{E}_1$, it is easy to check that 
	\begin{align*}
		\psi \mathcal{E}_t&=8\cdot D\left[8^{t-1} D^{t-1}(\psi)\right] =8\cdot D[\psi\cdot \mathcal{E}_{t-1}] \\
		&=[\mathcal{E}_{t-1}\cdot 8D(\psi)+\psi\cdot 8 D(\mathcal{E}_{t-1})] \\
		&=[\mathcal{E}_{t-1}\cdot \psi E_2+\psi\cdot 8D(\mathcal{E}_{t-1})].
	\end{align*}
	Dividing through by $\psi$, we arrive at the desired conclusion.
	\end{proof}

\noindent
The next result can be regarded as a generalization of the derivative rules \eqref{diffeq}.

\begin{theorem} We have the differential equation
	$$D(\mathbb{E}_{2t-2}(q))=t(2t+1)\cdot\mathbb{E}_{2t}(q)-3\mathbb{E}_2(q)\cdot\mathbb{E}_{2t-2}(q).$$
\end{theorem}
\begin{proof} For brevity, write $g(t):=\frac1{4^t(2t+1)!}$. Once more, we utilize 
	$8\frac{D(\psi)}{\psi}=E_2=\mathcal{E}_1=24\,\mathbb{E}_2$. From \cite[equation (3.3)]{AOS} we discern $\mathbb{E}_{2a}=g(a)D(\mathcal{E}_a$). 
Together with Lemma~\ref{diff-for-F}, these allow
	\begin{align*}
		D(\mathbb{E}_{2t-2})&=g(t-1)\cdot D(\mathcal{E}_{t-1})
		=\frac18\,g(t-1)[\mathcal{E}_t-E_2\cdot \mathcal{E}_{t-1}] \\
		&=\frac18\,g(t-1)\left[\frac1{g(t)}\,\mathbb{E}_{2t}-E_2\cdot \frac1{g(t-1)}\mathbb{E}_{2t-2}\right] \\
		&=\frac18[8t(2t+1)\cdot\mathbb{E}_{2t}-E_2\cdot \mathbb{E}_{2t-2}] \\
		&=\frac18[8t(2t+1)\,\mathbb{E}_{2t}-(24\,\mathbb{E}_2)\cdot \mathbb{E}_{2t-2}].
	\end{align*}
	Therefore, we gather $D(\mathbb{E}_{2t-2})=t(2t+1)\cdot\mathbb{E}_{2t}- 3\mathbb{E}_2\cdot\mathbb{E}_{2t-2}$. 
\end{proof}

\noindent
We chose to record the following result $\partial_{E_2}\mathcal{U}_{2t}=\sum_{j=1}^t\epsilon_j\mathcal{U}_{2t-2j}$ (similarly for $\mathcal{U}_{2t}^{\star}$) which determines 
$\mathcal{U}_{2t}$ up to an $E_2$ independent term (the pure modular part) and expressing it recursively in terms of those of lower weight members.

\begin{proposition} Preserving the notation and definitions from \eqref{MO_U} and \eqref{AAT_U},
we have that
$$\partial_{E_2}\mathcal{U}_{2t}(q)=-\frac1{12}\sum_{j=1}^t \frac{\mathcal{U}_{2t-2j}(q)}{j^2\binom{2j}j} \qquad \text{and} \qquad
\partial_{E_2}\mathcal{U}_{2t}^{\star}(q)=\frac1{12}\sum_{j=1}^t \frac{\mathcal{U}_{2t-2j}^{\star}(q)}{j^2\binom{2j}j}.$$
\end{proposition}

\begin{proof} Denote $\mathfrak{U}^{\star}(x):=\sum_{t\geq0}\mathcal{U}_{2t}^{\star}(q)\,x^t$. Direct calculation implies that
$$\log \mathfrak{U}^{\star}(x)=\sum_{r\geq1}H_r(q)\frac{x^r}r \qquad \text{and} \qquad H_r(q)=\sum_{k\geq1}\frac{q^{rk}}{(1-q^k)^{2r}}.$$
Let $\mathbf{S}_j(q)=\sum_{m\geq1}\frac{m^jq^m}{1-q^m}$. It is proven in \cite[Example 7.1]{AAT} that 

\begin{align} \label{umbral}
(2r-1)!H_r&=\mathbf{S}(\mathbf{S}^2-1^2)(\mathbf{S}^2-2^2)\cdots(\mathbf{S}^2-(r-1)^2)
\end{align}
with the \emph{umbral} convention (i.e. $\mathbf{S}^m$ understood as $\mathbf{S}_m$ after expansion).

\smallskip
\noindent
Operating with the derivative $\partial_{E_2}$ on the equation $\log \mathfrak{U}^{\star}(x)=\sum_{r\geq1}H_r(q)\frac{x^r}r$ leads to
$\sum_t x^t\partial_{E_2}\mathcal{U}_{2t}^{\star}=\mathfrak{U}^{\star}\sum_r\frac{x^r}r\partial_{E_2}H_r$ from which we obtain
\begin{align*}
\partial_{E_2}\mathcal{U}_{2t}^{\star}&=\sum_{j=1}^t \frac{\mathcal{U}_{2t-2j}^{\star}}j\cdot \partial_{E_2}H_j \\
&=\sum_{j=1}^t \frac{\mathcal{U}_{2t-2j}^{\star}}{j(2j-1)!}\cdot \partial_{E_2}\left[\mathbf{S}\prod_{s=1}^{j-1}(\mathbf{S}^2-s^2)\right] \\
&=\sum_{j=1}^t \frac{\mathcal{U}_{2t-2j}^{\star}}{j(2j-1)!} \cdot (-1)^j\cdot\frac{1^22^2\cdots(j-1)^2}{24}
\end{align*}
where the two facts $E_2=1-24\,\mathbf{S}_1$ and the identity in \eqref{umbral} have been utilized. We gather 
$$\partial_{E_2}\mathcal{U}_{2t}^{\star}=\frac1{12}\sum_{j=1}^t \frac{(-1)^j\,\mathcal{U}_{2t-2j}^{\star}}{j^2\binom{2j}j}.$$
Starting with $\mathfrak{U}(x):=\sum_{t\geq0} \mathcal{U}_{2t}(q)\,x^a$, it is evident that $-\log \mathfrak{U}(-x)=\sum_rH_r(q)\frac{x^r}r$. Then, applying an analogous argument as in above, one can show the remaining claim.
\end{proof}

\begin{remark} The relation depicted in \eqref{umbral} has a natural generalization
$$(\beta t-1)!\sum_{k\geq1}\frac{q^{tk}}{(1-q^k)^{\beta t}}=\prod_{s=1}^{\beta t-1}(\mathbf{S}-t+s).$$
\end{remark}

\section{Concluding remarks}

\noindent
Define the modular form $G_2(q):=2E_2(q^2)-E_2(2)$. We observe that $\mathbb{C}[E_2,\Theta_{0,1},\Theta_{1,1}]=\mathbb{C}[E_2,G_2,E_4]$, which is due to the conversion formulas for the modular forms $G_2$ and $E_4$:
$$G_2(q)= \Theta_{0,1}(q) \qquad \text{and} \qquad E_4(q)=\Theta_{0,1}^2(q)+ 6\Theta_{1,1}(q).$$ 
Thus all quasimodular forms over the congruence subgroup $\Gamma_0(2)$ form the ring $\mathbb{C}[E_2,G_2,E_4]$. 
Consequently, there exist some constants $\widetilde{c}_c(\alpha,\beta,\gamma)$ such that $\widetilde{c}_c(0,0,0)=1$ and 
\begin{align*}
\widetilde{c}_c(\alpha,\beta,\gamma)
&=\frac{2\alpha+4\beta+8\gamma-1}{24}\widetilde{c}_c(\alpha-1,\beta,\gamma)+\frac{8\gamma-8\beta+7}{24}\widetilde{c}_c(\alpha,\beta-1,\gamma) \\
&  - \frac{\alpha+1}{12} \widetilde{c}_c(\alpha+1,\beta,\gamma-1) + \frac{\beta+1}6 \widetilde{c}_c(\alpha,\beta+1,\gamma-1) \\
& - \frac{4(\gamma+1)}3 \widetilde{c}_c(\alpha,\beta-3,\gamma+1). 
\end{align*}
This allows one to write an alternative formulation of Theorem~\ref{theorem3}.

\begin{theorem} If $t$ is a positive integer, then we have that
\begin{align} \label{alt_thm3} \frac{D^t(\theta_4(q))}{\theta_4(q)} 
=\sum_{\substack{\alpha, \beta, \gamma\geq0 \\  \alpha+\beta+2\gamma=t}}
\widetilde{c}_c(\alpha,\beta,\gamma)\cdot E_2^{\alpha}(q)\, G_2^{\beta}(q)\, E_4^{\gamma}(q).
\end{align}
\end{theorem}

\begin{proof} The proof is again by induction on $t$. We require a few basic calculations: 
\begin{align} \label{useful_1}
\frac{D(\theta_4)}{\theta_4}&=\frac{E_2-G_2}{24}, \qquad D(E_2)=\frac{E_2^2-E_4}{12},  \\
D(G_2)&=\frac{E_2G_2-2G_2^2+E_4}6, \qquad D(E_4)=\frac{E_2E_4-4G_2^3+3G_2E_4}3. \nonumber
\end{align}
The base case $t=1$ is covered in \eqref{useful_1} with $\widetilde{c}_c(1,0,0)=\frac1{24}$ and 
$\widetilde{c}_c(0,1,0)=-\frac1{24}$.
Assume \eqref{alt_thm3} holds true. The next steps are very similar to the other proofs in this paper, so these are omitted to avoid unduly replications.
\end{proof}

\smallskip
\noindent
We record the below results (with a great deal of overlaps with \cite{ALL} and \cite{Fu-Sellers}) which reveal that if we add infinite families of MacMahon's quasimodular forms $\mathcal{U}_{2t}(q),\ \mathcal{C}_{2t}(q)$ and $\mathcal{U}^{\star}_{2t}(q)$ then the outcomes are weight $0$ modular forms although apparently they are of weights $\leq 2t$.

\begin{proposition} \label{prop_in5} Adopt the notation $(a;q)_{\infty}=(1-a)(1-aq)(1-aq^2)\cdots$ and the short-hand $(q)_{\infty}=(q;q)_{\infty}$. We have the identities
\begin{align*}
\sum_{t\geq0} \mathcal{U}_{2t}(q) 
&=\frac{(q^6)_{\infty}}{(q)_{\infty}(q^2)_{\infty}(q^3)_{\infty}},  \qquad
\sum_{t\geq0} \mathcal{C}_{2t}(q) 
= \frac{(q^4)_{\infty}(q^6)^2_{\infty}}{(q)_{\infty}(q^3)_{\infty}(q^{12})_{\infty}}, \\
\sum_{t\geq0} \mathcal{U}^{\star}_{2t}(q) 
& =\sum_{n\geq1} \frac{(-1)^{n-1}(1-q^n)(1-q^{2n})q^{\binom{n}2}}{1-3q^n+q^{2n}}.
\end{align*}
\end{proposition}

\begin{proof}
From \cite[Corollary 2]{Andrews-Rose}, we obtain
\begin{align*}
\sum_{t\geq0} \mathcal{U}_{2t}(q) &= \frac1{(q)^3_{\infty}} \sum_{n\geq0} q^{\binom{n+1}2} \left((-1)^n(2n+1)\sum_{t\geq0}  \frac{(-1)^t\binom{n+t}{2t}}{2t+1}\right)
=\frac{\sum_{n\geq0}q^{\binom{n+1}2} - 3\sum_{n\geq0}q^{\binom{3n+2}2}}{(q)^3_{\infty}}.
\end{align*}
On the other hand, we can simplify
\begin{align*}
\frac{(q^6)_{\infty}}{(q)_{\infty}(q^2)_{\infty}(q^3)_{\infty}}
&=\frac1{(q)^3_{\infty}} \prod_{k\geq1} \frac{1-q^{6k}}{(1+q^k)(1+q^k+q^{2k})} =\frac1{(q)^3_{\infty}} \prod_{k\geq1} (1-q^k+q^{2k})(1-q^k). 
\end{align*}
It suffices to recognize the identity $\prod_{k} (1-q^k+q^{2k})(1-q^k)=\sum_{n}q^{\binom{n+1}2} - 3\sum_{n}q^{\binom{3n+2}2}$.

\smallskip
\noindent
Again, from \cite[Corollary 2]{Andrews-Rose}, we obtain
\begin{align*}
\theta_4(q)  \sum_{t\geq0} \mathcal{C}_{2t}(q) &= \theta_4(q) + \sum_{n\geq0} q^{n^2} \left((-1)^n(2n)\sum_{t\geq1}  \frac{(-1)^t\binom{n+t}{2t}}{n+t}\right) 
=\sum_{\mathbb{Z}}q^{(3n)^2}-\sum_{\mathbb{Z}} q^{(3n+1)^2}
\end{align*}
to compare against ($\theta_4(q)=\frac{(q;q)_{\infty}}{(-q;q)_{\infty}}$ is used here)
\begin{align*}
\frac{\theta_4(q)\,(q^4)_{\infty}(q^6)^2_{\infty}}{(q)_{\infty}(q^3)_{\infty}(q^{12})_{\infty}}
&=\frac{(q)_{\infty}(-q^3,q^3)_{\infty}}{1-q^{2k}+q^{4k}} =\sum_{\mathbb{Z}}q^{(3n)^2}-\sum_{\mathbb{Z}} q^{(3n+1)^2}.
\end{align*}
The conclusion is clear. The last identity follows directly from \cite[Proposition 4.1]{AAT} where
$$\mathcal{U}^{\star}_{2t}(q)=\sum_{n\geq1} \frac{(-1)^{n-1}(1+q^n)q^{\binom{n}2+tn}}{(1-q^n)^{2t}}.  \qquad \qquad \qquad \qedhere$$
\end{proof}

\smallskip
\noindent
The following is an immediate consequence of Proposition \ref{prop_in5}.

\begin{corollary}
If we let $\sum_{t\geq0} \mathcal{U}_{2t}(q)=\sum_{n\geq0} u(n)q^n$ and $\sum_{t\geq0} \mathcal{C}_{2t}(q)=\sum_{n\geq0} \kappa(n)q^n$, then
\begin{align*}
u(3n+2)\equiv 0 \pmod 3 \qquad \text{and} \qquad \kappa(9n+6)\equiv 0 \pmod 3.
\end{align*}
\end{corollary}

\noindent
In line with the above, consider the infinite series that is defined in \cite[eq'n (1)]{AAT2} for which the authors have found a single-sum representation as
$$\mathcal{U}_{2t}(2;q)=\frac{(q)_{\infty}}{(q^2)^2_{\infty}}\sum_{n\geq0}\binom{n+t}{2t} q^{\binom{n+1}2}.$$
We now state and prove our claim which may be regarded as a result of independent interest.

\begin{theorem} \label{Fibo} If we let $\sum_{n\geq0} y(n)\, q^n:=\sum_{t\geq0} \mathcal{U}_{2t}(2;q)$, then
$$y(n) \equiv (-1)^n \cdot \# \{ (r,s) \in \mathbb{Z}^2 :\, n = r^2 + s^2 \} \pmod 5.$$
\end{theorem}

\begin{proof} Using the basic facts $\sum_{t\geq0} \binom{n+t}{2t}=F_{2n+1}$ and $F_{2n+1}\equiv (-1)^n(2n+1) \pmod 5$,
where $F_k$ stands for the Fibonacci number, we proceed to compute modulo $5$:
\begin{align*}
\sum_{t\geq0}\mathcal{U}_{2t}(2;q) &= \frac{(q)_{\infty}}{(q^2)^2_{\infty}} \sum_{n\geq0}F_{2n+1}q^{\binom{n+1}2}
\equiv \frac{(q)_{\infty}}{(q^2)^2_{\infty}} \sum_{n\geq0}(-1)^n(2n+1)q^{\binom{n+1}2}
=\frac{(q)_{\infty}(q)^3_{\infty}}{(q^2)^2_{\infty}} = \theta^2_4(q).
\end{align*}
However, $\theta^2_4(q)=\sum_{r,s\in\mathbb{Z}}(-1)^{r+s} q^{r^2+s^2}=\sum_{n\geq0} [(-1)^n\cdot \#\{(r,s)\in\mathbb{Z}: n=r^2+s^2\}]q^n$.
\end{proof}

\smallskip
\noindent
We close this section and our paper with a remark (the simple proof is left for the reader) regarding the sequence $c_v(\alpha,\beta,\gamma)$ defined in equation \eqref{coefficientsV} and attributed to Ramanujan's $q$-series $V_{2t}(q)$. Namely, these coefficients carry a closed form when varying $\alpha$ while $\beta$ and $\gamma$ remain fixed. The merit or value of our formula can be viewed as a way of generating coefficients, in the expansion portrayed in Theorem \ref{theorem2}, of quasimodular components on the basis of the corresponding modular forms.

\begin{proposition} \label{vary_alpha} Letting $(x)_n=x(x+1)\cdots(x+n-1)$ be the Pochhammer symbol, it holds true that
$$\frac{c_v(\alpha,\beta,\gamma)}{c_v(0,\beta,\gamma)}=\frac{(1+4\beta+6\gamma)_{2\alpha}}{2^{\alpha}\,\alpha!}.$$
\end{proposition}

\end{document}